\documentclass[12pt]{amsart}

\usepackage{amsfonts,amsmath,amsthm}
\usepackage{latexsym}

\newtheorem{theorem}{Theorem}[section]

\newtheorem{remark}[theorem]{Remark}
\newtheorem{problem}[theorem]{Problem}

\newtheorem{conjecture}[theorem]{Conjecture}

\newtheorem{definition}{Definition}[section]

\theoremstyle{definition}
\theoremstyle{problem}
\theoremstyle{conjecture}

\newcommand{\dst}{\displaystyle}

\newcommand{\Co}{\ensuremath{\mathbb{C}}}

\def \C {\mathbb{C}}

\newcommand{\ac}{\ensuremath{\mathcal{A}}}
\newcommand{\bc}{\ensuremath{\mathcal{B}}}

\newcommand{\eb}{\ensuremath{\mathbf{e}}}
\newcommand{\fb}{\ensuremath{\mathbf{f}}}

\def \< {\langle}
\def \> {\rangle}

\newcommand{\ent}[1]{{\left[{#1}\right]}}
\newcommand{\abs}[1]{{\left|{#1}\right|}}
\newcommand{\scal}[1]{{\left\langle{#1}\right\rangle}}

\begin{document}

\title[The MUB problem in dimension 6]{The problem of mutually unbiased bases in dimension 6}

\author[P. Jaming]{Philippe Jaming}
\address{P.J.: Universit\'e d'Orl\'eans\\
Facult\'e des Sciences\\
MAPMO - F\'ed\'eration Denis Poisson\\ BP 6759\\ F 45067 Orl\'eans Cedex 2\\
France} \email{philippe.jaming@univ-orleans.fr}

\author[M. Matolcsi]{M\'at\'e Matolcsi}
\address{M. M.: Alfr\'ed R\'enyi Institute of Mathematics,
Hungarian Academy of Sciences POB 127 H-1364 Budapest, Hungary
Tel: (+361) 483-8302, Fax: (+361) 483-8333}
\email{matomate@renyi.hu}

\thanks{M. Matolcsi was supported by OTKA Grant No. K77748}

\author[P. M\'ora]{P\'eter M\'ora}
\address{P. M. (and M. M. part time): BME Department of Analysis,
Egry J. u. 1, H-1111 Budapest, Hungary}
\email{morapeter@gmail.com}


\begin{abstract}
We outline a discretization approach to determine the maximal
number of mutually unbiased bases in dimension 6. We describe the
basic ideas and introduce the most important definitions to tackle
this famous open problem which has been open for the last 10
years. Some preliminary results are also listed.
\end{abstract}

\maketitle

{\it Dedicated to Prof. Warwick de Launey
on the occasion of
 his 50th birthday}

\bigskip

\bigskip


{\bf Keywords and phrases.} {\it  Mutually unbiased bases, complex
Hadamard matrices}

\section{Introduction}

This paper is based on the talk given by the second author at the
International Conference on Design Theory and Applications, NUI,
Galway, July 1-3, 2009.

\medskip

The notion of mutually unbiased bases (MUBs) constitutes a basic
concept of Quantum Information Theory and plays an essential role
in quantum-tomography \cite{Iv,WF}, quantum cryptography
\cite{BPT,BB,Re}, the mean king problem \cite{AE} as well as in
constructions of teleportation and dense coding schemes
\cite{wer}.

\medskip

Recall that two orthonormal bases of $\Co^d$,
$\ac=\{\eb_1,\ldots,\eb_d\}$ and $\bc=\{\fb_1,\ldots,\fb_d\}$ are
said to be \emph{unbiased} if, for every $1\leq j,k\leq d$,
$\abs{\scal{\eb_j,\fb_k}}=\dst\frac{1}{\sqrt{d}}$. A set
$\bc_0,\ldots\bc_m$ of orthonormal bases is said to be
\emph{(pairwise) mutually unbiased} if every two of them are
unbiased. It is well-known ( see e.g. \cite{BBRV,BBELTZ,WF}) that
the number of mutually unbiased bases (MUBs) in $\Co^d$ cannot exceed
$d+1$. It is also known that $d+1$ such bases can be constructed
if the dimension $d$ is a prime or a prime power ({\it see e.g.}
\cite{BBRV,Com0,Com1,Com2,Iv,KR,WF}). Such a set of $d+1$ MUBs in dimension $d$ is called a {\it complete
set}. If the dimension $d=p_1^{\alpha_1}\dots p_k^{\alpha_k}$ is composite then very little
is known except for the fact that there are at least $p_j^{\alpha_j}+1$ mutually
unbiased bases in $\C^d$ where $p_j^{\alpha_j}$ is the smallest of the prime-power divisors.
Thus, the first case where the largest number of mutually
unbiased bases is unknown is $d=6$:

\begin{problem}\label{MUB6problem}\ \\
What is the maximal number of pairwise mutually unbiased bases in
$\Co^6$?
\end{problem}

Although this famous open problem has received considerable
attention over the past few years
(\cite{BBELTZ,config,ujbrit,arxiv,Msz,Skinner}), it remains wide open.
Since $6=2\times 3$, we know that there are at least $3$ mutually
unbiased bases in $\C^6$ (see also \cite{Za,arxiv} for {\it infinite families} of MUB-triplets), but so far tentative numerical evidence
\cite{config,ujbrit,numerical,Za} suggests that there are no more
than 3, a fact apparently first conjectured by Zauner \cite{Za}.

\begin{conjecture}\label{MUB6conj}\ \\
The maximal number of pairwise mutually unbiased bases in $\Co^6$
is 3.
\end{conjecture}

One reason for the slow progress is that mutually unbiased bases
are naturally related to \emph{complex Hadamard matrices}.
Indeed, if the bases $\bc_0,\ldots,\bc_m$ are mutually
unbiased we may identify each
$\bc_l=\{\eb_1^{(l)},\ldots,\eb_d^{(l)}\}$ with the \emph{unitary}
matrix
$$
[H_l]_{k,j}=\ent{\scal{\eb_k^{(l)},\eb_j^{(0)}}_{1\leq k,j\leq
d}},
$$
{\it i.e.} the $k$-th row of $H_l$ consists of the
coordinates of the $k$-th vector of $\bc_l$ in the basis $\bc_0$.
(Throughout the paper the scalar product $\scal{.,.}$ of $\C^d$ is
linear in the first variable and conjugate-linear in the second.
Note also that for convenience of computer programming we use the
unconventional definition that the {\it rows} of the matrices
correspond to the vectors of the bases.) With this convention,
$H_0=I$ the identity matrix and all other matrices are unitary
and have entries of modulus $1/\sqrt{d}$. Therefore, the matrices
$\sqrt{d}H_l$ have all entries of modulus 1 and complex orthogonal
rows (and columns). Such matrices are called \emph{complex
Hadamard matrices}. It is clear that the existence of a family of
mutually unbiased bases $\bc_0,\ldots,\bc_m$ is thus equivalent to
the existence of a family of complex Hadamard matrices
$\sqrt{d}H_1,\ldots, \sqrt{d}H_m$ such that for all
 $1\leq j\not=k\leq
m$, $\sqrt{d}H_jH_k^*$ is again a complex Hadamard matrix. In such
a case we will say that these complex Hadamard matrices are {\it
mutually unbiased}.

\medskip

In particular, the existence of 4 MUBs $\bc_0, \bc_1, \bc_2, \bc_3$ in dimension 6 is equivalent to the existence of 3 mutually unbiased Hadamard matrices $H_1, H_2, H_3.$ Therefore, in an attempt to prove that no collection of 4 MUBs exist in dimension 6 it is enough to prove that no triplet of mutually unbiased Hadamard matrices exist. This will be the core of our argument in this note.

\medskip

A complete classification of complex Hadamard matrices, however, is only
available up to dimension 5 ({\it see} \cite{haagerup})
which allows for a complete classification of MUBs ({\it see} \cite{BWB}). 
The
classification in dimension 6 is still out of reach despite recent
efforts \cite{BN,Msz,Skinner}. This is one of the reasons for
Problem \ref{MUB6problem} to be difficult.

\medskip

In this paper we outline a discretization approach that is likely
to lead to the proof of Conjecture \ref{MUB6conj} in the near
future. Once all the ideas are properly implemented in a computer
code, an {\it exhaustive} search will be carried out to prove
Conjecture \ref{MUB6conj}. We will include all the basic
definitions and ideas as well as some preliminary results.
We remark here that a partial result of this approach has already been completed:
in \cite{arxiv} we assumed that the first Hadamard matrix $H_1$ comes from the two-parameter Fourier family of complex Hadamard matrices,
and we proved by discretization and an exhaustive computer search that in such a case a MUB-quartet $I, H_1, H_2, H_3$ cannot exist. In this paper, however, we tackle the general case, so that we must consider $H_1$ as {\it any} complex Hadamard matrix of dimension 6. This complicates matters quite considerably as the number of cases to check after the discretization increases by orders of magnitude.
As an optimistic note let us recall here that the non-existence of a projective plane of
order 10 was also proved by an exhaustive computer search
\cite{lam}.

\section{Discretization}\label{sec2}

The proof proceeds by contradiction, via a discretization scheme.
Assume that there exists a collection of 4 MUB's in $\Co^6$.
Equivalently, there exist $6\times 6$ complex Hadamard matrices
$A,B,C$ having all entries of modulus $1$, such that the rows (and
thus the columns) are complex orthogonal, and we have the unbiased
condition: for any two rows $u,v$ coming from different matrices
we have $|\langle u, v \rangle |=\sqrt{6}$. (Recall that for the
purposes of this note the {\it rows} of the matrices correspond to
the vectors of the bases.) In such a case the orthonormal bases
$\frac{1}{\sqrt{6}}A, \frac{1}{\sqrt{6}}B, \frac{1}{\sqrt{6}}C$
accompanied with the identity matrix $Id$ correspond to a family
of 4 MUB's. We assume that such matrices $A, B, C$ exist and try
to reach a contradiction.

\medskip

After multiplying rows and columns by appropriate scalars if
necessary, we can assume that all coordinates of the first row and
column of $A$ are 1's, and all coordinates of the first column of
all other matrices are 1's (i.e. we assume that all
vectors in the bases $A,B,C$ have first coordinate 1, and the
first vector in basis $A$ is an all 1's vector). All the other
coordinates in the matrices are complex numbers of modulus 1, i.e.
they are of the form $e^{2\pi i \rho}$ with $\rho\in [0,1)$. We
will use a discretization approach. Let $N$ be a positive integer,
called the {\it discretization parameter}. We partition the
interval $[0,1)$ into $N$ sub-intervals
$I_0^{(N)},I_1^{(N)},\ldots, I_{N-1}^{(N)}$ of equal length, {\it
i.e.} $I_j^{(N)}=[j/N,(j+1)/N)$. (Other partitions are also
possible, but this seems most convenient for programming.) Now,
any entry $e^{2\pi i \rho}$ in any of the matrices $A, B, C$ will
be represented by the integer $j$ if $\rho \in I_j^{(N)}$ (note
that $0\le j\le N-1$). This means: whenever we see an entry $j$
somewhere in a matrix then we conclude that the original phase
$\rho$ must lie somewhere in the interval $I_j^{(N)}$. We have
{\it no more and no less information than this}. We also agree
that the first coordinate of each row will be represented by $0$,
keeping in mind that it represents exactly 1, without error (and
not the interval $I_0^{(N)}$).

\medskip

In short: we will exclusively be dealing with row vectors of the
form
\begin{equation}\label{basicform}u=(0, j_1, j_2, j_3, j_4, j_5)
\end{equation}
 where $0\le j_k\le N-1$ and
the first coordinate 0 represents 1 without error, while the other
coordinates $j_k$ mean that the actual entry $\rho_k$ falls into
the interval $I_{j_k}^{(N)}$. In notation, the original matrix
will be denoted by $A$, while its representative integer matrix
will be denoted by $\tilde A$. The entries of $A$ will be denoted
by $\rho_{m,k}$, while those of $\tilde A$ will be denoted by
$j_{m,k}$.

\medskip

There are altogether $N^5$ vectors of the form \eqref{basicform}.

\medskip

Also, there is a natural ordering among these vectors: $u\le v$ if
and only if it is so in lexicographical order. We will use this
ordering throughout this note.

\section{The search for the discretized Hadamard matrix $\tilde A$}\label{sec:A}

The matrix $\tilde A$ is an integer matrix with first row and
column consisting of 0's and the core of the matrix containing
integers between 0 and $N-1$. We introduce the following
definition:

\begin{definition}
Given an integer matrix $\tilde A$ with first row and column
consisting of 0's and the core of the matrix containing integers
$j_{m,k}$ between 0 and $N-1$, we will say that $\tilde A$ is an
{\em $N$-discretized representative of a complex Hadamard matrix}
if there exists a complex Hadamard matrix $A$ with entries
$\rho_{m,k}$ such that $\rho_{m,k}\in I_{j_{m,k}}$. In notation
$\tilde A \in HAD_N$, where $HAD_N$ denotes the set of
$N$-discretized representatives of complex Hadamard matrices.
\end{definition}

\medskip

The aim of this section is to describe an algorithm to efficiently
search for all possible matrices $\tilde A \in HAD_N$. Upon strong
numerical evidence \cite{Skinner}, it is conjectured that the
manifold of $6\times 6$ complex Hadamard matrices is
4-dimensional. Therefore we expect that the cardinality of $HAD_N$
will be approximately $c N^4$ for some constant $c$. Nevertheless,
the task of finding all possible $\tilde A$ is daunting at first
glance. There are $N^{25}$ possible $N$-discretized matrices
altogether, and we must select the ones belonging to $HAD_N$. The
number $N^{25}$ is of course astronomical even for $N\approx 50$,
but we will see that with an intelligent approach the task can
still be carried out.

\medskip

There are a few properties we can assume about $\tilde A$ without
loss of generality. We already assumed that the first row and
column consist of 0's. We can also assume that both the rows and
columns are arranged so that they increase with respect to
lexicographical order. This can be arranged by repeated
permutation of rows and columns. (This is not entirely trivial
because ordering the rows lexicographically can actually spoil
such an ordering of the columns and vice versa. However, if one
writes out the matrix entries row-after-row in one 36-long row
vector, then it is clear that this vector will decrease
lexicographically irrespectively of whether you make an ordering
of the rows or the columns. Therefore such a repeated ordering of
rows and columns will terminate in finite steps, and will produce
a matrix such that both the rows and the columns increase in
lexicographical order.) This automatically implies that the
entries of the second row and second column are both monotonically
increasing. This is a very convenient property, because it
restricts the possibilities for the second row and column quite
strongly.

\medskip

We can also assume that the second row is less than or equal to
the second column in lexicographical order (this can be arranged
by transposition of $\tilde A$ if necessary).

\medskip

We must make use of the fact that the rows (and columns) of $A$
are complex orthogonal to each other. The first row and column of
$\tilde A$ consist of 0's (representing the entry 1 in $A$,
without error). Therefore, we have 5 unknown rows and columns of
$\tilde A$. All of these rows and columns have the form
\eqref{basicform}. The orthogonality condition with the first row
(and column) makes it natural to introduce the following
definition:

\begin{definition}\label{ortn}
We will say that a vector $u$ of the form \eqref{basicform}
belongs to $ORT_{N}$ if there exist $\phi_k\in I_{j_k}$ such that
$1+\sum_{k=1}^5 e^{2i\pi \phi_k}=0$.
\end{definition}

Note that $ORT_N$ is a ``small'' subset of all the vectors of form
\eqref{basicform}, containing only those vectors which represent
vectors being orthogonal to the vector $(1,1,1,1,1,1)$. Clearly,
all rows and columns of $\tilde A$ must belong to $ORT_N$.
Therefore it is very important to determine the set $ORT_N$ as
precisely as we can. We achieve this by the following ``check the
descendants'' method.

\medskip

Let $u=(0, j_1, j_2, j_3, j_4, j_5)$, and let $r_{j_k}$ denote the
midpoint of the interval $I_{j_k}$ (the superscript $N$ has been
dropped from the notation for convenience). If $u\in ORT_N$ then
the trivial error bound ({\it see} Lemma 3.1 in \cite{arxiv}) gives
\begin{equation}\label{trivialest}
|1+\sum_{k=1}^5 e^{2\pi i r_{j_k}}| \le \frac{5\pi}{N}.
\end{equation}

This is too crude, but we can iterate it to the ``children'' of $u$.
Namely, assume that the numbers $\phi_k$ exist as in Definition
\ref{ortn}. For each interval $I_{j_k}$ the value $\phi_k$ must
lie in either the left or the right half of $I_{j_k}$. There are
32 choices, according to whether we consider the left or the right
half of each interval $I_{j_k}$. These choices are called the 32
``children'' of $u$. Clearly, at least one of these children needs to
satisfy \eqref{trivialest} with $\frac{5\pi}{2N}$ on the right
hand side (and its own midpoints substituted to the left hand
side, of course, instead of $r_{j_k}$). If none of the children
satisfy this, then $u$ can be discarded. Of course we iterate this
to grandchildren, and so on, down to 7-8 generations. The vector
$u$ survives this test if it has at least one surviving descendant
in each generation.

\begin{remark}
{\rm The set $ORT_N$ is clearly invariant under permutations of
the last 5 coordinates $j_1, j_2, j_3, j_4, j_5$. Therefore it
makes sense to introduce the set $ORT_{N, mon}$ of vectors in
$ORT_N$ with monotonically increasing coordinates. To save time,
in the actual computer code we first find the vectors of $ORT_{N,
mon}$ by the method above, and then we permute the last 5
coordinates to arrive at the set $ORT_N$.}
\end{remark}

\begin{remark}
{\rm There exists also an improved error bound ({\it see} Lemma 3.2 in
\cite{arxiv}). It is somewhat slower to check by computer and it
is reasonable to believe that we arrive at the same set $ORT_N$ by
applying either error bounds.}
\end{remark}

\begin{remark}
{\rm We have implemented a computer code for selecting the set
$ORT_N$. For example, for $N=17$ we have $|ORT_N|=58450$, for
$N=19$, $|ORT_N|=82630$, and for $N=53$, $|ORT_N|=1875110$.
Experience shows that the set $ORT_N$ is unexpectedly large if $N$
is divisible by 2 or 3. Therefore, we have mainly restricted our
attention to $N$ being a prime.}
\end{remark}

\begin{remark}
{\rm The optimal choice of $N$ seems to be crucial for the success
of the project. Clearly, if $N$ is too small then the error bounds
are not good enough and we will not reach a contradiction in the
forthcoming steps ({\it see} Section \ref{sec3} below). However, if $N$
is too large then the size of the sets $ORT_N$ and correspondingly
$HAD_N$ will be far too large to be manageable. At present we
believe that the optimal choice of $N$ is around $N\approx 50$.}
\end{remark}

\medskip

Let us turn back to the construction of $\tilde A$. All rows and
columns must come from $ORT_N$, and they must be pairwise
$N-orthogonal$ in the following sense:

\begin{definition}
We will say that the vectors $u=(0, j_1, j_2, j_3, j_4, j_5)$ and
$v=(0, m_1, m_2, m_3, m_4, m_5)$ are {\em $N$-orthogonal} if there
exist numbers $\phi_k$ and $\psi_k$ in the intervals $I_{j_k}$ and
$I_{m_k}$, such that $1+\sum_{k=1}^5 e^{2i\pi (\phi_k-\psi_k)}=0$.
\end{definition}

This property is clearly shift-invariant in the sense that it only
depends on the values $(j_1-m_1, \dots j_5-m_5)$ modulo $N$. We
can therefore take $m_1=\dots=m_5=0$ and correspondingly $v_0=(0,
0 \ldots, 0)$, (where the last 5 coordinates represent the
interval $I_0$, of course) and define the set $ORT_{eps, N}$ as
the set of vectors of the form \eqref{basicform} which are
$N$-orthogonal to $v_0$. (The notation $ORT_{eps, N}$ indicates
that the vector $v_0$ contains an ``epsilon'' of error, because the
last 5 coordinates represent the interval $I_0$ and not the exact
number 1.) With this notation the shift-invariance means that $u$
and $v$ will be $N$-orthogonal if and only if the vector
$(j_1-m_1, \dots j_5-m_5) (mod \ N)$ is in $ORT_{eps, N}$.

\medskip

Having constructed the set $ORT_{N}$ previously, it is now easy to
obtain $ORT_{eps, N}$. Indeed, by definition a vector
$u=(0,j_1,\dots j_5)$ can only be $N$-orthogonal to $v_0$ if there
exist numbers $\phi_k$ in the intervals $I_{j_k}$ and $\psi_k$ in
$[0,\frac{1}{N})$, such that $1+\sum_{k=1}^5 e^{2i\pi
(\phi_k-\psi_k)}=0$. But then the numbers $\phi_k-\psi_k$ must
fall in the intervals $I_{j_k-\epsilon_k}$ where $\epsilon_k$ is
either 0 or 1, and hence the vector $u_\epsilon=(0,
j_1-\epsilon_1,\ldots,j_5-\epsilon_5)$ is in $ORT_{N}$.

\medskip

Therefore, $ORT_{eps, N}$ will consist of all the vectors of the
form $u^\epsilon=(0,j_1+\epsilon_1,\ldots,j_5+\epsilon_5)$, where
$\epsilon_k$ is 0 or 1, and the vector $(0,j_1,\ldots,j_5)$ is in
$ORT_{N}.$

\medskip

\begin{remark}
{\rm Each $u\in ORT_N$ gives rise to 32 different $u^\epsilon$
above. One could therefore expect that the size of $ORT_{eps, N}$
will be nearly 32 times the size of $ORT_N$. This is not so,
however, because there will be many coincidences. Experience shows
that the size of $ORT_{eps, N}$ is approximately 4 times the size
of $ORT_N$, regardless of the value of $N$.}
\end{remark}

\medskip

Now we are ready to conduct a search for the possible matrices
$\tilde A$. The first row and column are full of 0's. We then
build up the matrix with a row-by-column approach. We fit in the
second row, then the second column, then the third row, then the
third column, etc. At each step we must consider that:

\medskip

-- each row and column must come from $ORT_N$.

\medskip

-- each row (resp. column) must be lexicographically larger than
any previous rows (resp. columns). In particular, the entries of
the second row and column are monotonically increasing, i.e. they
belong to $ORT_{N, mon}$.

\medskip

-- the second column must be lexicographically larger than or
equal to the second row.

\medskip

-- each row (resp. column) must be $N$-orthogonal to any previous
rows (resp. columns). This is equivalent to the fact that the
pairwise differences of the rows (resp. columns) modulo $N$ must
be contained in $ORT_{eps, N}$.

\medskip

-- each row (resp. column) must be compatible with the already
existing entries of the matrix ({\it e.g.} when we fit in the fourth
row, then its first 3 coordinates are already fixed because the
first three columns of the matrix have already been filled out
previously).

\medskip

We have implemented a computer code which executes the search as
described above. The running time is still reasonable, within 1-4
days, depending on $N$. However, the number of selected matrices
$\tilde A$ is unexpectedly large. It is in the range of $10^9 - 5
\cdot 10^{10}$ as $N$ ranges from $17$ to $53$. Let $PREHAD_N$
denote the set of matrices obtained by this search. Clearly,
$HAD_N \subset PREHAD_N$.

\medskip

Have we made all possible restrictions so as to list {\it
exclusively} the matrices $\tilde A$ belonging to $HAD_N$? In
other words, is it true that $HAD_N = PREHAD_N$? It turns out that
this is not the case, and there is an important possibility for
further pruning. Consider a matrix $\tilde A \in PREHAD_N$. There
are 25 non-trivial entries in $\tilde A$ (the first row and column
being trivial), all of which represent intervals $I_{j_{m,k}}$ of
length $1/N$. Once again we can ``check the descendants'' of $\tilde
A$. That is, we can take left or right halves of each 25 intervals
$I_{j_{m,k}}$, and therefore consider the $2^{25}$ children of
$\tilde A$. Obviously, at least one of these children need to
satisfy {\it stricter pairwise orthogonality conditions} of rows
and columns. If none of the children do, then $\tilde A$ can be
discarded, i.e. it does not belong to $HAD_N$. Of course, checking
$2^{25}$ children is very slow, but if one proceeds row-by-row
then only a few thousand children need to be actually checked. We
have not rigorously implemented this step in our computer code.
Nevertheless, preliminary results suggest that only a small
fraction of the matrices in $PREHAD_N$ will pass this test, i.e.
$HAD_N$ will be significantly smaller in size than $PREHAD_N$.
This is very important for the running time of the overall
algorithm, as the size of $HAD_N$ should definitely be kept in the
range $10^8 - 10^9$ even for $N\approx 50$.

\section{Stage 2: vectors unbiased to $\tilde A,$ and reaching a
contradiction}\label{sec3}

Let us fix a matrix $\tilde A\in HAD_N$. We want to prove that the
pair $(Id, \tilde A)$ cannot be extended by matrices $\tilde B,
\tilde C$ so as to meet all orthogonality and unbiasedness
conditions. The rows of $\tilde B$ and $\tilde C$ are of the form
\eqref{basicform} and they must be ``unbiased'' to all six rows of
$\tilde A$. Therefore, as a next step, we must obtain a list of
all such vectors.

\begin{remark}
{\rm We are actually free to use a different discretization
parameter $N'$ for the matrices $\tilde B$ and $\tilde C$. It may
well reduce the running time if we use optimal choices for $N$ and
$N'$. Experience shows ({\it see} \cite{arxiv}) that it makes sense to
choose $N'$ considerably smaller than $N$. However, for the sake
of simplicity we will keep $N=N'$ throughout this note.}
\end{remark}

As the first row of $\tilde A$ is invariably $(0,0,0,0,0,0)$
(representing 1's in the first row of $A$, without error) it makes
sense to introduce the following definition:

\begin{definition}
We will say that a vector $u=(0,j_1, j_2, j_3, j_4, j_5)$ belongs
to the set $UB_{N}$ if there exist $\phi_k\in I_{j_k}$ such that
$|1+\sum_{k=1}^5 e^{2i\pi \phi_k}|=\sqrt{6}$. We will say that $u$
belongs to $UB_{N, mon}$ if the coordinates of $u$ are
monotonically increasing.
\end{definition}

\medskip

The set $UB_N$ can be constructed in a similar way as $ORT_N$.
With $r_{j_k}$ denoting the midpoint of the interval $I_{j_k}$ the
trivial estimate gives
\begin{equation} \left | \ |1+\sum_{k=1}^5 e^{2i\pi
r_{j_k}}|-\sqrt{6} \right |\le \frac{5\pi}{N}. \end{equation}
This is too crude, of course, and the descendants of $u$ need to
be checked for some 7-8 generations.

\begin{remark}
{\rm Once again, the set $UB_N$ is invariant under the permutation
of the last 5 coordinates $j_1, j_2, j_3, j_4, j_5$. Therefore, in
practice, we first check monotonically increasing vectors only,
and obtain $UB_{N,mon}$. Then we permute the coordinates to obtain
$UB_N$.}
\end{remark}

\begin{remark}
{\rm The set $UB_N$ is much larger than $ORT_N$. This can be
expected because orthogonality of complex vectors induces two
conditions (the real part and imaginary part both being zero)
while unbiasedness only induces one condition.}
\end{remark}

\begin{remark}
{\rm We have implemented a code for listing the set of vectors
$UB_N$. For example, for $N=17$ we have $|UB_N|=479340$, while for
$N=19$, $|UB_N|=764060$.}
\end{remark}

We will also need a set $UB_{eps, N}$ which is analogous to
$ORT_{eps, N}$.

\begin{definition}
We will say that the vectors $u=(0, j_1, j_2, j_3, j_4, j_5)$ and
$v=(0, m_1, m_2, m_3, m_4, m_5)$ are {\em $N$-unbiased} if there
exist numbers $\phi_k$ and $\psi_k$ in the intervals $I_{j_k}$ and
$I_{m_k}$, such that $|1+\sum_{k=1}^5 e^{2i\pi
(\phi_k-\psi_k)}|=\sqrt{6}$.
\end{definition}

This property is again shift-invariant in the sense that it only
depends on the values $(j_1-m_1, \dots j_5-m_5)$ modulo $N$. We
can therefore take $m_1=\dots=m_5=0$ and correspondingly $v_0=(0,
0 \ldots, 0)$, (where the last 5 coordinates represent the
interval $I_0$, of course) and define the set $UB_{eps, N}$ as the
set of vectors of the form \eqref{basicform} which are
$N$-unbiased to $v_0$.  With this notation the shift-invariance
means that $u$ and $v$ will be $N$-unbiased if and only if the
vector $(j_1-m_1, \dots j_5-m_5) (mod \ N)$ is in $UB_{eps, N}$.

\medskip

Let the vector $v$ of the form \eqref{basicform} be any row of
either $\tilde B$ or $\tilde C$. As the first row of $\tilde A$ is
invariably $(0,0,0,0,0,0)$ (representing the 1's in the first row
of $A$, without error), we conclude that $v\in UB_N$. Let $a_2,
\dots a_6$ denote the last five rows of $\tilde A$. Then, by
definition, the differences $v-a_j$ modulo $N$ must belong to
$UB_{eps, N}$ for all $j=2, \dots 6$. Let $UB_{\tilde A}$ denote
the set of vectors $v$ which satisfy these conditions. The
notation $UB_{\tilde A}$ reflects that these are the vectors which
are ``unbiased'' to all rows of $\tilde A$. By what has been said
above, all rows of $\tilde B$ and $\tilde C$ must belong to
$UB_{\tilde A}$.

\begin{remark}
{\rm We have implemented a code to obtain the set $UB_{\tilde A}$.
Experience shows that the size of $UB_{\tilde A}$ is largely
independent of the choice of $\tilde A$, and it has between $10^3 -
10^4$ vectors as $N$ ranges from 17 to 53.}
\end{remark}

Having constructed $UB_{\tilde A}$ we must show that $\tilde B$
and $\tilde C$ cannot be built from these vectors satisfying all
orthogonality and unbiased conditions.

\medskip

Consider the vectors in $UB_{\tilde A}$ and try to build the
matrix $\tilde B$ out of them. This means that we need to find 6
vectors $b_1,\dots b_6$ such that the pairwise differences
$b_k-b_m$ modulo $N$ all belong to $ORT_{eps, N}$. Counting
constraints and parameters one would expect that only a finite
number of triplets of MUB's $(Id, A, B)$ exists, and a MUB-pair
$(Id, A)$ can generically {\it not} be extended to a triplet $(Id,
A, B)$. This would give us hope that a contradiction is reached
most of the times while trying to build $\tilde B$. However,
recent results \cite{arxiv, ujbrit} show that infinite families of
MUB-triplets do exist. Numerical practice also shows that the
matrix $\tilde B$ can indeed be built from the vectors of
$UB_{\tilde A}$ for all $\tilde A$. Therefore, we do not get an
immediate contradiction. Instead, for each $\tilde B$ we must go
on and select the vectors $UB_{\tilde A, \tilde B}$ which are
unbiased to all rows of $\tilde A$ and $\tilde B$, and we must try
to build a matrix $\tilde C$ out of the vectors $UB_{\tilde A,
\tilde B}$. The contradiction is reached only at this point. That
is, if $N$ is large enough the matrix $\tilde C$ cannot be
constructed from $UB_{\tilde A, \tilde B}$ to meet all
orthogonality conditions. Experience shows that $N$ must be larger
than 30 to reach a contradiction. This part of the project is
currently under implementation. It would be desirable to reach a
contradiction for each $\tilde A$ within a few seconds of
computing time.

\medskip

For the overall success of the project two tasks need to be
considered in the near future. One is the implementation of the
ideas of the last paragraph of Section \ref{sec:A} to bring down
the number of possible $\tilde A$'s to the region $10^8-10^9$. The
other is to reach a contradiction for each $\tilde A$ within a few
seconds of computing time by not being able to construct $\tilde
B$ and $\tilde C$. A nice feature of the overall project is that
once the algorithm is completed, it is very easy to distribute the
calculations among several hundreds of computers, and thus
reducing the running time by 2-3 orders of magnitude.

\medskip

Finally, we remark that the entire discretization procedure
described above has already been completed in \cite{arxiv} in the
restricted setting when $A$ is assumed to belong to the Fourier
family $F(a,b)$ of complex Hadamard matrices.

\begin{theorem}\label{no4s}{\rm [Theorem 1.4 in
\cite{arxiv}]}    \ \\
None of the pairs $\bigl(Id, F(a,b)\bigr)$ of mutually unbiased
orthonormal bases can be extended to a quartet $\bigl(Id,
F(a,b),B,C\bigr)$ of mutually unbiased orthonormal bases.
\end{theorem}

In that case we used the discretization parameters $N=180$ for
$\tilde A$ and $N'=19$ for $\tilde B$ and $\tilde C$. Due to some
well-known equivalence relations only a few hundred possible
discretized matrices $\tilde A$ needed to be considered, and a
contradiction was quickly reached for all of them. The
documentation of that search is available at \cite{web}. The
difficulty in the general case is that the number of matrices
$\tilde A$ becomes very large if $N$ is chosen large, while if $N$
is small then a contradiction is reached very slowly (or not
reached at all!) in the second stage of the search.

\end{document}